\documentclass[11pt]{article}

\usepackage{rotating}
  \usepackage{graphicx}
  \usepackage{psfrag}
  \usepackage{amsfonts}
\usepackage{amsmath,amssymb}
\usepackage{url}

\usepackage{natbib}
\usepackage{color}

\newtheorem{theorem}{Theorem}[section]
\newtheorem{problem}{Problem}[section]
\newtheorem{property}{Property}[section]
\newtheorem{lemma}[theorem]{Lemma}
\newtheorem{proposition}[theorem]{Proposition}
\newtheorem{corollary}[theorem]{Corollary}
\newtheorem{definition}[theorem]{Definition}
\newtheorem{remark}[theorem]{Remark}

\newenvironment{example}[1][Example]{\begin{trivlist}
\item[\hskip \labelsep {\bfseries #1.}]}{\end{trivlist}}

\newcommand{\qed}{\nobreak \ifvmode \relax \else
      \ifdim\lastskip<1.5em \hskip-\lastskip
      \hskip1.5em plus0em minus0.5em \fi \nobreak
      \vrule height0.75em width0.5em depth0.25em\fi}

\def\comment#1{\textit{[#1]}}
\def\comment#1{}

\begin{document}

\title{Combinatorics of least squares trees}
\author{Radu Mihaescu\thanks{Department of Mathematics, UC Berkeley} \and Lior Pachter\footnotemark[1] \footnote{lpachter@math.berkeley.edu}}
\date{\today}
\maketitle

\begin{abstract}

A recurring theme in the least squares approach to phylogenetics has been the discovery of elegant combinatorial formulas for the least squares estimates of edge lengths. These formulas have proved useful for the development of efficient algorithms, and have also been important for understanding connections among popular phylogeny algorithms. For example, the selection criterion of the neighbor-joining algorithm is now understood in terms of the combinatorial formulas of Pauplin for estimating tree length.

We highlight a phylogenetically desirable property that weighted least squares methods should satisfy, and provide a complete characterization of methods that satisfy the property. The necessary and sufficient condition is a multiplicative four point condition that the the variance matrix needs to satisfy.
The proof is based on the observation that the Lagrange multipliers in the proof of the Gauss--Markov theorem are tree-additive. Our results generalize and complete previous work on ordinary least squares, balanced minimum evolution and the taxon weighted variance model. They also
provide a time optimal algorithm for computation.
\end{abstract}

\section{Introduction}

The least squares approach to phylogenetics was first suggested by
Cavalli-Sforza \& Edwards \cite{Cavalli-Sforza1967} and Fitch \& Margoliash \cite{Fitch1967}. The  precise problem formulated in \cite{Cavalli-Sforza1967} was Problem \ref{prob:ols}:

\begin{definition}[Pair-edge incidence matrix]
Given a phylogenetic $X$-tree $T$ with edge set $E$ and $|X|=n$ (see \cite{Semple2003} for basic definitions), the {\em pair-edge incidence matrix} of $T$ is the ${n
  \choose 2} \times |E|$  matrix
  \[ (S_T)_{ij,e} = \left\{ \begin{array}{l} 1 \mbox { \normalsize if } e \in E \mbox{ \normalsize is an edge on the
    path between } i\mbox{ \normalsize and } j,\\0 \mbox { \normalsize otherwise.} \end{array} \right. \]
\end{definition}
\begin{definition}[Tree-additive map]
Let $T$ be a phylogenetic $X$-tree. A dissimilarity map $D$ is $T$-additive if for some vector $l \in \mathbf{R}^{|E|}$,
\begin{equation}
D_{ij} = (S_T  l)_{ij}.
\end{equation}
\end{definition}

\begin{problem}[Ordinary least squares]
\label{prob:ols}
Find the phylogenetic $X$-tree $T$ and $T$-additive map $\hat{D}$ that minimizes
\begin{equation}
\sum_{i,j \in {X \choose 2}} (D_{ij}-\hat{D}_{ij})^2.
\end{equation}
\end{problem}
For a fixed tree, the solution of Problem \ref{prob:ols} is a linear algebra problem (Theorem \ref{thm:least_squares_soln}). However Rzhetsky \& Nei \cite{Rzhetsky1993} showed that the {\bf O}rdinary {\bf L}east {\bf S}quares edge lengths could instead be computed using elegant and efficient combinatorial formulas. Their result was based on an observation of Vach \cite{Vach1989}, namely that OLS edge lengths obey the desirable {\em {\bf I}ndependence of {\bf I}rrelevant {\bf P}airs} property (our choice of terminology is inspired by social choice theory \cite{Ray1973}):
\begin{property}[IIP]

Let $T$ be a phylogenetic $X$-tree and $e$ an edge in $T$. A linear edge length estimator for $e$ is a linear function from dissimilarity maps to the real numbers, i.e. $\hat{l}_e=\sum_{ij}p_{ij}D_{ij}$. We say that such an estimator satisfies  the IIP property if $p_{ij}=0$ when the path from $i$ to $j$ in $T$ (denoted $\overline{i,j}$ ) does not contain either of $e$'s endpoints.
\end{property}

In other words, the IIP property is equivalent to the statement that the sufficient statistic for the least squares estimator of the length of $e$ is a projection of the dissimilarity map onto the coordinates
given by pairs of leaves whose joining path contains at least one endpoint of $e$. It has been shown that this crucial property is satisfied not only by ordinary least squares (OLS) estimators, but also by specific instances of {\em {\bf W}eighted {\bf L}east {\bf S}quares} estimators (e.g., \cite{Semple2003b}).

\begin{problem}[Weighted least squares]
\label{prob:wls}
Let $T$ be a phylogenetic $X$-tree and $D$ be a dissimilarity map. Find
the $T$-additive map $\hat{D}$ that minimizes
\begin{equation}
\label{eq:WLS}
\sum_{i,j \in {X \choose 2}} \frac{1}{V_{ij}}\left( D_{ij}-\hat{D}_{ij}\right)^2.
\end{equation}
\end{problem}

The {\em variance matrix} for weighted least squares is the ${n \choose 2} \times {n \choose 2}$ diagonal
matrix $V$ whose diagonal entries are the $V_{ij}$. Note that $V$ can also be regarded as a dissimilarity map and we will do so in this paper. Weighted least squares for trees was first suggested in \cite{Fitch1967} and \cite{Hartigan1967}, with the former proposing specifically $V_{ij}=D_{ij}^2$.
\begin{theorem}[Least squares solution]
\label{thm:least_squares_soln}
The solution to Problem \ref{prob:wls} is given by
$\hat{D}=S_T  \hat{l}$ where
\begin{equation}
\label{eq:least_squares}
\hat{l} = (S_T^t   V^{-1}  S_T)^{-1}  S_T^t  V^{-1}
 D.
\end{equation}
\end{theorem}
We note that The OLS problem reduces to the case $V=I$. The statistical significance of the variance matrix together with a statistical interpretation of Theorem
\ref{thm:least_squares_soln} is provided in Section 2.

It follows from (\ref{eq:least_squares}) that the lengths of the edges
in a weighted least squares tree are linear combinations of the entries of the
dissimilarity map. A natural question is therefore which variances matrices $V$ result in edge length estimators that satisfy the IIP property? Our main result is an answer to this question in the form of a characterization (Theorem \ref{thm:main}): a WLS model is IIP if and only if the variance matrix is semi-multiplicative. We show that such matrices are good approximations to the variances resulting from popular distance estimation procedures. Moreover, we provide combinatorial formulas that describe the WLS edge lengths under semi-multiplicative variances (Equation \ref{eq:WLSInternal}), and show that they lead to optimal algorithms for computing the lengths (Theorem \ref{thm:compute}).

The key idea that leads to our results is a connection between Lagrange multipliers arising in the proof of the Gauss--Markov theorem and the weak fundamental theorem of phylogenetics that provides a combinatorial characterization of tree-additive maps (Remark \ref{rem:main_point}). This explains many isolated results in the literature on least squares in phylogenetics; in fact, as we show in the section "The multiplicative model and other corollaries", almost all the known theorems and algorithms about least squares estimates of edge lengths follow from our results.
\section{BLUE Trees}\label{sec:BLUE}
The foundation of least squares theory in statistics is the Gauss--Markov
theorem. This theorem states that the {\bf B}est {\bf L}inear {\bf U}nbiased {\bf E}stimator for a linear combination of the edge lengths, when the errors have
zero expectation, is a least
squares estimator. We explain this theorem in the context of Problem
\ref{prob:wls}.

\begin{lemma}
For any phylogenetic $X$-tree $T$, the matrix $S_T$ is full rank.
\end{lemma}
{\bf Proof}: We show that for any $e \in E$, the vector
$f_e=(0,\ldots,1,\ldots,0)$ of size $|E|$ with a $1$ in the $e$-th position
and $0$ elsewhere lies in the row span of $S$. Choose any $i,j,k,l \in
X$ such that the paths from $i$ to $j$ and from $k$ to $l$ do not
intersect, and the intersection of the paths from $i$ to $j$ and from
$k$ to $l$ is exactly the edge $e$. Note that
\begin{equation}
\frac{1}{2}\sum_e \left( S_{ik,e}+S_{jl,e}-S_{ij,e}-S_{kl,e}\right) = f_e.
\end{equation}

\begin{theorem}[Gauss--Markov Theorem]
\label{thm:Gauss}
Suppose that $D$ is a random dissimilarity map of the form $D=S_T  l+\epsilon$
where $T$ is a tree, and $\epsilon$ is a vector of random
variables satisfying $E(\epsilon)=0$ and $Var(\epsilon)=V$ where
$V$ is an invertible variance-covariance matrix for $\epsilon$.

Let $M(S_T^t)$ be the linear space generated by the columns of $S_T^t$ and $f \in M(S_T^t)$. Then $f^t  \hat{l} = p^t  D$ (where $\hat{l}$ given
by (\ref{eq:least_squares})) has minimum variance
among the linear unbiased estimators of $f^t  l$.

\end{theorem}
{\bf Proof}:
Observe that the problem of finding $p$ is equivalent to solving a
constrained optimization problem:
\begin{equation}
\mbox{min } p^t  V  p \, \mbox{ subject to } \, S_T^t  p = f.
\end{equation}
The first condition specifies that the goal is to minimize the variance;
the second constraint encodes the requirement that the estimator is unbiased.
Using Lagrange multipliers, it is easy to see that the minimum variance
unbiased estimator of $f^t  l$ is the unique
vector $p$ satisfying
\begin{eqnarray}
\label{eq:lagrange}
V  p & = & S_T  \mu \mbox{ for some } \mu \in {\bf R}^{|E|}, \\
\label{eq:lagrange1}
S_T^t  p & = & f.
\end{eqnarray}
In other words
\begin{eqnarray}
\begin{pmatrix}
V & -S_T\\
 S^t_T & 0
\end{pmatrix}\begin{pmatrix} p \\ \mu \end{pmatrix} & =
& \begin{pmatrix} 0\\f \end{pmatrix}
\end{eqnarray}
\begin{eqnarray}
\Rightarrow \, \begin{pmatrix} p\\ \mu \end{pmatrix} & = &
\begin{pmatrix}
V^{-1}S_TU^{-1}S^t_TV^{-1} & (U^{-1}S^t_TV^{-1})^t\nonumber\\
-U^{-1}S^t_TV^{-1} & U^{-1}
 \end{pmatrix}
\begin{pmatrix}
  0\\ f \end{pmatrix}
\end{eqnarray}

where $U=S^t_TV^{-1}S_T$.

The Gauss--Markov Theorem can also be proved directly using linear
algebra, but the Lagrange multiplier proof has two advantages: First, it
provides a description of $p$ different from
(\ref{eq:least_squares}) that is simpler and more informative. Secondly,
the technique is general and can be
used in many similar settings to find minimum variance unbiased
estimators. Hayes and Haslett \cite{Hayes1999} provide pedagogical
arguments in favor
of Lagrange multipliers for interpreting least squares
coefficients and discuss the origins of this approach in applied statistics \cite{Matheron1962}.

In phylogenetics, Theorem \ref{thm:Gauss} (and its proof) are useful because for each edge $e$, the vector $f_e$
in the standard basis for $M(S_T^t)$ is associated with a vector $p$ such
that $p^tD$
is the best linear unbiased estimator for the length of $e$. Similarly, the tree length is estimated
from $f_T=(1,1,\ldots,1)$ which is also in $M(S_T^t)$.
Condition (\ref{eq:lagrange}) is particularly interesting because it says that there exists some $T$-additive map $\Lambda= S^{t}_T \mu=V  p$, whose (possibly negative) edge lengths are given by the Lagrange multipliers $\mu$. 

\begin{figure}[!h]
\begin{center}
\includegraphics[scale=0.47]{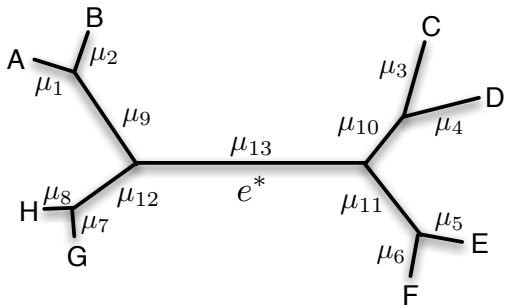}
\end{center}
\caption{The Lagrange tree $\Lambda$ for an IIP weighted least squares estimator for the central edge $e^*$ of a complete binary tree with 8 leaves.
In Proposition, \ref{prop:8tax} $X=\{A,B,C,D,E,F,G,H\}$, whereas in the proof of Theorem \ref{thm:main} the leaf labels represent clades. The IIP property means that the WLS estimate $\hat{l}_{e^*}$ does not depend on $D_{AB},D_{CD},D_{EF}$ or $D_{GH}$.}
\label{fig:8taxa}
\end{figure}
The following theorem provides a combinatorial characterization of tree-additive maps, and hence of the {\em Lagrange tree} $\Lambda$:
\begin{definition}[Weak four point condition]
A dissimilarity map $D$ satisfies the weak four point condition if for
any $i,j,k,l \in X$, two of the following three linear forms are equal:
\begin{equation}
D_{ij}+D_{kl}, \quad D_{ik}+D_{jl}, \quad D_{il}+D_{jk}.
\end{equation}
\end{definition}
\begin{theorem}[Weak fundamental theorem of phylogenetics]
\label{thm:four_point}
A dissimilarity map $D$ is tree-additive if and only if it satisfies the weak four point condition.
\end{theorem}
Theorem \ref{thm:four_point} was first proved in \cite{Patrinos1972}. For a recent exposition see Corollary 7.6.8 of \cite{Semple2003} where it is derived using the theory of group-valued dissimilarity maps.
We note that the pair of equal quantities in the four point condition define the topology of a quartet. Furthermore the topology of the tree is defined uniquely by the topologies of all its quartets. We again refer the reader to \cite{Semple2003} for details.

The {\em Lagrange equations} (\ref{eq:lagrange}) and (\ref{eq:lagrange1}) together with Theorem \ref{thm:four_point} form the mathematical basis for our results:
\begin{remark}
\label{rem:main_point}
Condition (\ref{eq:lagrange}) specifies that $Vp$ must be a $T$-additive map. It follows
that $Vp$
satisfies the weak four point condition. In other words, (\ref{eq:lagrange})
amounts to a combinatorial characterization of $Vp$, and hence $p$.
Condition (\ref{eq:lagrange1}) imposes a
normalization requirement on $p$. Together these conditions are useful
for finding $p$, and also for understanding its combinatorial properties.
\end{remark}

The structure of the Lagrange tree in the case of $OLS$ is the middle quartet of the tree shown in Figure 1. It immediately reveals interesting properties of the estimator. For example the fact that it is a tree on four taxa implies the IIP property. The content of \cite[Appendix 2]{Desper2004} is that for tree length estimation under the balanced minimum evolution model, the Lagrange tree is the star tree. In fact, we will see that most of the known combinatorial
results about least squares estimates of edge and tree lengths can be
explained by Remark \ref{rem:main_point} and interpreted in terms of the structure of the Lagrange tree.

\section{Main Theorem}\label{sec:model}

Our main result is a characterization of IIP WLS estimators. In the sections that follow we will see that the IIP property for WLS is not only biologically desirable, but also statistically motivated and algorithmically convenient.
We begin by introducing some notation and concepts that are necessary for stating our main theorem.
\begin{definition}[Clade]
A clade of a phylogenetic $X$-tree $T$ is a subset $A \subset X$ such that there exists
an edge in $T$ whose removal induces the partition $\{A,X \setminus A \}$.
We also use clade to mean the induced topology $T|_A$.
\end{definition}

Given a dissimilarity map $D$ and a variance matrix $V$, we set
\begin{eqnarray*}
D_{AB} & := &\sum_{a\in A, b\in B}V_{ab}^{-1}D_{ab}, \mbox{ and}\\
Z_{AB} &:= &\sum_{a\in A, b\in B}V_{ab}^{-1}.
\end{eqnarray*}
where $A,B$ are disjoint clades. If $e_1, \ldots, e_k \in E(T)$  form a path with ends determining clades $A$ and $B$, then by the notation $D_{e_1 \cdots e_k}$ and $Z_{e_1\cdots e_k}$ we mean $D_{AB}$ and $Z_{AB}$ respectively.
Note that if $e$ is an edge in a tree $T$ then (\ref{eq:lagrange},\ref{eq:lagrange1}) imply that the Lagrange tree for any WLS estimate of $e$ satisfies $\Lambda_e = f_e$.
\begin{definition}[Semi-multiplicative map]
A dissimilarity map $D$ is semi-multiplicative with respect to disjoint clades $A,B$
if for any $a_1,a_2\in A$ and $b_1,b_2\in B$
\begin{equation}
\label{eq:mult}
D_{a_1b_1}D_{a_2b_2}=D_{a_1b_2}D_{a_2b_1}.
\end{equation}
We say that $D$ is semi-multiplicative with respect to $T$ if for any pair of disjoint clades $A,B$, not defined by the same edge of $T$, (\ref{eq:mult}) holds.
\end{definition}
\begin{lemma}
$D$ is semi-multiplicative if and only if every clade
$A$ of $T$ has the property that for any $A' \subset A$, and any clade $B$ disjoint from $A$ and induced by a different edge, for all $x \in B$,
\begin{equation}
\label{eq:star}
Z_{\{x\}A'}/Z_{\{x\}A}=\xi^B_{A'A},\end{equation}
where $\xi^B_{A'A}$ does not depend on $x$.
\end{lemma}
It is an easy exercise to prove that $A$ satisfies (\ref{eq:star}) for all relevant $B$ if and only if (\ref{eq:star}) holds for the the two clades disjoint from $A$ and defined by the two edges adjacent to the edge defining $A$.

The semi-multiplicative condition is slightly weaker than $\log D$ being tree-additive. Indeed, removing the requirement that the clades $A,B$ are defined by different edges of $T$ leaves one one with a multiplicative analog of the four-point condition. By Theorem \ref{thm:four_point}, this is equivalent to $D_{ij}=\prod_{e\in \overline{i,j}}w(e)^{-1}$ for some $w:E(T) \rightarrow \mathbf{R}_{+}$  \cite{Gill2008}.
\begin{theorem}[Characterization of IIP WLS estimators]
\label{thm:main}
A WLS edge length estimator for an edge in a tree $T$ has the IIP property if and only if the variance matrix is semi-multiplicative with respect to $T$.
\end{theorem}
The proof of the theorem reduces to the WLS solution for the length of an edge in a tree with at most eight leaves (edge $e^*$ in Figure \ref{fig:8taxa}):
\begin{proposition}
\label{prop:8tax}
Let $T$ be the phylogenetic $X$-tree shown in Figure \ref{fig:8taxa}.
The Lagrange tree $\Lambda=S_T \mu$ for the WLS problem of estimating the length of the edge $e^*$ satisfies the property that $\mu_1=-\mu_2$, $\mu_3=-\mu_4$, $\mu_5=-\mu_6$ and $\mu_7=-\mu_8$. Furthermore, these Lagrange multipliers and the remaining ones $\mu_9,\ldots,\mu_{13}$ can be computed by solving $\mu = (S_T^t V^{-1}S_T)^{-1}f_{e^*}$.
\end{proposition}
{\bf Proof}:
Using the notation of Figure \ref{fig:8taxa}, with the convention that the edge labeled by $\mu_i$ is $e_i$, it follows from (\ref{eq:lagrange1}) that $\Lambda_{e_i}=0$ for $i=1,2,9$. But $\Lambda_{e_i}=\Lambda_{e_i e_j}+\Lambda_{e_i e_k}$ for $\{i,j,k\}=\{1,2,9\}$, which implies that $\Lambda_{e_i e_j}=0$ $\forall i,j\in \{1,2,9\}$. Therefore $V_{AB}^{-1}\Lambda_{AB}=V_{AB}^{-1}(\mu_{1}+\mu_{2})=0$ and the result follows. The arguments for $e_3,e_4$, $e_5,e_6$ and $e_7,e_8$ are identical. The complete solution for the $\mu$ for a given $V$ is given by $\mu = (S_T^t V^{-1}S)^{-1}f_{e^*}$, which reduces to the inversion of a $13 \times 13$ matrix. 

Note that the proof only uses the fact that $e_1,e_2$ are adjacent leaf edges not adjacent to $e^*$. The conclusion $\mu_{e_1}=-\mu_{e_2}$ will hold identically in any tree for a pair of edges of this type.\\
\\
{\bf Proof of Theorem \ref{thm:main}}:
We begin by showing that if $V$ is semi-multiplicative then the WLS edge length estimators have the IIP property. This calculation involves showing that for any phylogenetic $X$-tree $T$ and edge $e^* \in T$, the Lagrange tree for $e^*$ is the tree in Figure \ref{fig:8taxa}, where $A$,$B$,$C,D,E,F,G,H$ are clades with the property that their intra-clade Lagrange multipliers are zero.

Let $e_1,\ldots ,e_k$, with $k\leq 8$, be the edges of $T$ such that either $d(e^*,e_i)=2$ or $d(e^*,e_i)<2$ and $e_i$ is a leaf edge. For $i\in \{1,\ldots,k\}$, let $C_i$ be the clade defined by $e_i$ such that $e^*\not\in C_i$. Let $T^{/e^*}$ to be the phylogenetic $X^{/e^*}$-tree, where $X^{/e^*}=\{C_1,\ldots,C_k\}$, with topology induced by $T$ in the natural way (see Figure \ref{fig:8taxa}). Set $V^{/e^*}$ be the diagonal variance matrix on pairs of nodes in $X^{/e^*}$ given by $V^{/e^*}_{C_i C_j}=Z_{C_iC_j}^{-1}$.

If $\mu^{/e^*}$ are the Lagrange multipliers and $\Lambda^{/e^*}$ is the Lagrange tree given by estimating $\hat{l}_{e^*}$ for topology $T^{/e^*}$ and variance $V^{/e^*}$ then the $T$-additive map given by $\Lambda=S_T^t \mu$ with
\begin{equation}\mu_e=\begin{cases} \mu^{/e^*}_e \text{ if } e\in E(T^{/e^*}),\\ 0 \text{ otherwise.}\end{cases}\end{equation} satisfies the Lagrange equations for $T$. Thus $\mu$ are the Lagrange multipliers for $\hat{l}_{e^*}$ and $\hat{l}_{e^*}=\Lambda^t  V^{-1}  D$.

We let $\Lambda^{/e*}$, $Z^{/e*}$ denote the natural correspondents of $\Lambda$ and $Z$ for the problem of estimating $\hat{l}_{e^*}$ from and $V^{/e^*}$ and $T^{/e^*}$. It is an easy exercise to check that for all $e\in E(T^{/e^*})$, we have $Z^{/e*}_e=Z_e$ and $\Lambda^{/e*}_e=\Lambda_e$. This implies that $\Lambda_e=f_e$ for all $e\in E(T^{/e^*})$, i.e. the Lagrange equation (\ref{eq:lagrange1}) is satisfied for $e\in E(T^{/e^*})$.

Now consider edge $e\in C_1$. We need to verify that $\Lambda_e=0$. Since $\Lambda_{ij}=0$ for all $i,j\in C_1$, $\Lambda_e=\Lambda_{e \cdots e_2}+\Lambda_{e \cdots e_9}$. Now for all $i\in C_1$ and $j\in C_2$, $\Lambda_{ij}=\mu_1+\mu_2=0$, so $\Lambda_{e \cdots e_2}=0$. Finally let $A'\subset A$ be the clade defined by $e$ and let $A''$ be the clade defined by $e_9$ which does not intersect $A$. The fact that $V$ is semi-multiplicative implies that for any taxon $x\in A''$
\begin{equation}
Z_{\{x\}A'}/Z_{\{x\}A}=\xi^{C_1}_{A'A}
\end{equation} 
where $\xi_{A'A}$ does not depend on the taxon $x$. This implies $\Lambda_{e \cdots e_9}=\xi^{C_1}_{A'A}\Lambda_{e_1 \cdots e_9}=0$ by the proof of Proposition \ref{prop:8tax}.

Since $\mu_e=0$ for all $e\not\in T^{/e^*}$, it is enough to show that $\Lambda^{/e^*}$ satisfies the IIP property. This follows from Proposition \ref{prop:8tax}. Therefore, $V$ has the IIP property with respect to $T$, i.e. $\Lambda_{ij}=0$ for all $i,j\in X$ such that $\overline{i,j}$ does not intersect $e^*$.

This concludes the proof for the "if" part of Theorem \ref{thm:main}. For the "only if" direction, we will prove by induction that (\ref{eq:star}) is satisfied by all clades $A$ of $T$, and thus the variance $V$ is semi-multiplicative with respect to $T$. The base case is provided by clades formed by a single leaf, for which (\ref{eq:star}) holds vacuously. 

\begin{figure}[!h]
\begin{center}
\includegraphics[scale=0.47]{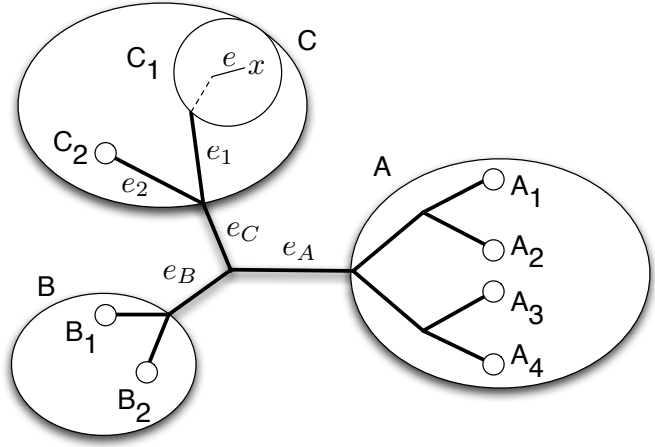}
\end{center}
\caption{Configuration of the induction in the proof that IIP WLS models are semi-multiplicative.}

\label{fig3}
\end{figure}

For the induction step, suppose clades $A$ and $B$ both satisfy (\ref{eq:star}), and that they are defined by adjacent edges $e_A$ and $e_B$ (see Figure 2). Let $e_C$ be the other edge adjacent to $e_A$ and $e_B$ and let $C=X\setminus (A\cup B)$ be the clade it defines. We would like to prove that the clade $(A\cup B)$ also satisfies (\ref{eq:star}). If $|C|=1$, this holds vacuously. We may therefore assume that there exist two more edges $e_1,e_2$ incident with $e_C$. Let $C_i\subset C$ be the clade defined by $e_i$, for $i=1,2$. It suffices to prove that $(A\cup B)$ satisfies (\ref{eq:star}) with respect to $C_1$ and $C_2$. Notice that $A$ and $B$ already satisfy (\ref{eq:star}) with respect to $C_1$ and $C_2$. Therefore it is enough to show that
\begin{equation}\frac{Z_{\{x\}A}}{Z_{\{x\}(A\cup B)}}=\xi^{C_1}_{A (A\cup B)}\end{equation}
is the same for all $x\in C_1$, and similarly for all $x\in C_2$.

Now consider the problem of estimating $\hat{l}_{e_A}$. Let $\mu$ be the corresponding Lagrange multipliers and $\Lambda=S_{T}\mu$ be the Lagrange tree they define. By the IIP property, $\Lambda$ defines an identically zero tree additive map on the clade $C$. Therefore the edge lengths corresponding to this map are all zero. This implies $\mu_e=0$ for all $e\in E(C), e\neq e_1,e_2$, and also $\mu_{e_1}+\mu_{e_2}=0$.

Let $A_1,\ldots, A_k$, with $k\leq 4$ and $B_1,\ldots, B_t$, with $t\leq 2$, be the sub-clades of $A$, respectively $B$, corresponding to nodes of $T^{/e_A}$. Then for any $x\in C_1$ and $y\in A_i$, and $z\in B_j$, $\Lambda_{xy}=\Lambda^{/e_A}_{C_1 A_i}$ does not depend on $x,y$ and $\Lambda_{xz}=\Lambda^{/e_A}_{C_1 B_j}$ does not depend on $x,z$.

Now pick $x\in C_1$ and let $e$ be the leaf edge adjacent to it. Then $\Lambda_e=0$. Since all Lagrange multipliers are 0 inside the clade $C_1$, $\Lambda_e=\Lambda_{e\ldots e_1}=\Lambda_{e\ldots e_2}+\Lambda_{e\ldots e_c}$. Since $\mu_{e_1}+\mu_{e_2}=0$, $\Lambda_{e\ldots e_2}=0$. Thus $\Lambda_{e\ldots e_C}=\Lambda_{\{x\}A}+\Lambda_{\{x\}B}=0$.
Equivalently,
\begin{eqnarray}
\sum_{i=1}^k Z_{\{x\},A_{i}}\Lambda^{/e_A}_{C_1 A_i}+\sum_{j=1}^t Z_{\{x\},B_{j}}\Lambda^{/e_A}_{C_1 B_j}=0 \Leftrightarrow\nonumber\\
 Z_{\{x\},A}\sum_{i=1}^k \xi^{C_1}_{A_iA}\Lambda^{/e_A}_{C_1 A_i}+Z_{\{x\},B}\sum_{j=1}^t \xi^{C_1}_{B_jB}\Lambda^{/e_A}_{C_1 B_j}=0
\end{eqnarray}
This imposes a linear equation on $Z_{\{x\}A}$ and $Z_{\{x\}B}$ whose coefficients do not depend on $x$. Thus the following also does not depend on $x$:
\begin{equation}
\xi^{C_1}_{A(A\cup B)}=\frac{Z_{\{x\}A}}{Z_{\{x\}(A\cup B)}}=\frac{Z_{\{x\}A}}{Z_{\{x\}A}+Z_{\{x\}B}}.\end{equation}

\section{An optimal algorithm for WLS edge lengths}
\label{sec:algorithm}
\begin{theorem}[Computing WLS edge lengths]
\label{thm:compute}
Let $D$ be a dissimilarity map and $V$ an IIP variance matrix. The set of all WLS edge lengths estimates for a tree $T$ can be computed in $O(n^2)$ where $n$ is the number of leaves in $T$.
\end{theorem}

\begin{figure}[!h]
\begin{center}
\includegraphics[scale=0.47]{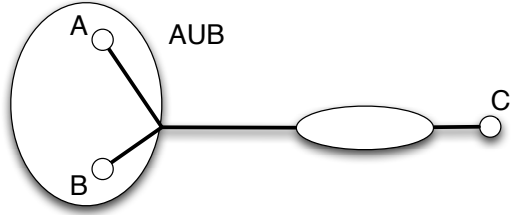}
\end{center}
\caption{Configuration of the dynamic programming recursion for computing WLS edge lengths. $A,B$ and $A \cup B$ are clades, and $C$ is a clade disjoint from $A \cup B$.The oval in the middle represents the rest of the tree.
}
\label{fig2}
\end{figure}

{\bf Proof}: It is apparent from the proof of Theorem \ref{thm:main} that all one needs in order to compute the WLS edge lengths are the values of $D_{AB}$ and $Z_{AB}$, where $A$ and $B$ are disjoint clades of $T$. We define the height of a tree to be the distance between its root and its farthest leaf, where the root is taken to be the closest endpoint of the edge defining the clade. Thus the height of a clade formed by just one leaf is $0$.

Now consider the configuration in Figure \ref{fig2}. The clades $A,B,C$ are all pairwise disjoint and $A$ and $B$ are adjacent. It is easy to see that $A\cup B$ form a clade for which \begin{eqnarray}
Z_{A\cup B,C}& = & Z_{AC}+Z_{BC},\\
D_{A\cup B,C} & = & (D_{AC}Z_{AC}+D_{BC}Z_{BC})/Z_{A\cup B,C}.
\end{eqnarray}

Therefore one needs only constant time to compute $D_{A\cup B,C}$ and $Z_{A\cup B,C}$ if $D_{AC}$,$Z_{AC}$,$D_{CB}$ and $Z_{CB}$
are known. Clearly, there are $O(n)$ clades since there are $O(n)$ edges, and thus there are $O(n^2)$ pairs of disjoint clades. We can compute $D_{AB}$ and $Z_{AB}$ for all pairs $AB$ through a simple dynamic program. We start with pairs of trees of height $0$, for which the values of $D$ and $Z$ are trivially given by $\delta$ and $V^{-1}$. After round $2t$ of the algorithm we will know $D_{AB}$ and $Z_{AB}$ for all disjoint pairs $A,B$ of height at most $t$ and after round $2t+1$ we know $D_{AB}$ and $Z_{AB}$ for all disjoint pairs $A,B$ of height $t+1$ and $t$ respectively. The algorithm clearly  requires constant time per clade pair. Subsequently, all $O(n)$ edge lengths can be computed in constant time per edge: the calculation of each edge length involves only a constant number of multiplications and one matrix inversion (of size at most $13 \times 13$). Thus the algorithm is optimal since its running time is proportional to the size of the input.

We note that many algorithms have been proposed for computing WLS edge lengths for certain specific models (these are discussed in the next section). Existing approaches rely on different recursive schemes that lead to markedly different algorithms. Some attempt to reduce the size of the problem by agglomerating leaves (\cite{Denis2003}); others start with a star topology and gradually extend it by refining internal nodes (\cite{Vach1989}). In fact, all these methods implicitly compute Lagrange multipliers in a recursive way, and dealing directly with Lagrange multipliers may in many cases clarify the exposition and suggest simplified implementations. As we can see from the above theorem however, once one has the closed form expressions for the edge lengths, these inductive arguments can be easily replaced by our dynamic program.

\section{The multiplicative model and other corollaries}\label{sec:corollaries}
In this section we begin by giving formulas for the WLS edge lengths assuming a a {\em tree-multiplicative} variance matrix, i.e. $V_{ij}=\prod_{e\in \overline{i,j}}w_e^{-1}$ for some $w:E(T) \rightarrow \mathbf{R}_{+}$. Throughout the section, $e^* \in E(T)$ denotes the edge for which the WLS length is being computed. If $e^*$ is an internal edge then $A,B,C,D$ are the adjacent clades. In the case that $e^*$ is adjacent to a leaf, that leaf is labeled $i$ and the adjacent clades $A,B$.

\begin{proposition}\label{lem:WLS_edges}
If $V$ is a tree-multiplicative variance matrix then the WLS edge length of an internal edge is
\begin{eqnarray}\label{eq:WLSInternal}
2\hat{l}_{e^*} & = & \frac{Z_{AD}+Z_{CB}}{Z_{A\cup B,C\cup D}}(D_{AC}+D_{BD})\nonumber\\
& + & \frac{Z_{AC}+Z_{DB}}{Z_{A\cup B,C\cup D}}(D_{AD}+D_{BC})\nonumber\\
& - & D_{AB}-D_{CD}.
\end{eqnarray}
If $e^*$ is adjacent to a leaf then the WLS length is 
\begin{equation}\label{eq:WLSLeaf}
2\hat{l}_{e^*}=D_{Ai}+D_{Bi}-D_{AB}.
\end{equation}
\end{proposition}

At first glance these formulas may seem surprising, but the derivation is straightforward after solving for the Lagrange multipliers.

{\bf Proof:} By the results of the previous section, it is enough to verify that the Lagrange equations hold. By Proposition \ref{prop:8tax} this is equivalent to verifying that the Lagrange equations hold for $T^{/e^*}$ and $V^{/e^*}$, which is a simple exercise left to the reader.

We now present a number of previous results about least squares that can be interpreted (and in some cases completed) using Theorems \ref{thm:main}, \ref{thm:compute}, and Lemma \ref{lem:WLS_edges}. All the models we discuss are special cases of the multiplicative variance model and all of our statements can be easily proven by substituting the appropriate form of $V$ into (\ref{eq:WLSInternal},\ref{eq:WLSLeaf}).

\begin{example}[Ordinary least squares]
\end{example}
This is the first model considered for least squares phylogenetics, and is the most studied model for edge and tree length estimation. It corresponds to the variance matrix equal to the identity matrix. 
\begin{corollary}[Rzhetsky \cite{Rzhetsky1993}]
The ordinary least squares estimate $p^tD = f_e^t(S_T^t  S_T)^{-1} S_T^tD$  for the
length of edge $e$ is given by
\begin{eqnarray}
2\hat{l}_{e^*} & = & \frac{n_A n_D+n_B n_C}{(n_A+n_B)(n_C+n_D)}(D_{AC}+D_{BD})\nonumber\\
& + & \frac{n_A n_C+n_B n_D}{(n_A+n_B)(n_C+n_D)}(D_{AD}+D_{BC})\nonumber\\
& - & D_{AB}-D_{CD},
\end{eqnarray}
where $n_A,n_B,n_C$ and $n_D$ are the number of leaves in the clades $A,B,C$ and $D$, and $D_{AC}=\sum_{a\in A,c\in C}D_{ac}$. If $e^*$ is a leaf edge, $\hat{l}_{e}$ is given by:
\begin{equation}\label{lem:WLSLeaf1}
2\hat{l}_{e^*}=D_{Ai}+D_{Bi}-D_{AB}.
\end{equation}
\end{corollary}
Our algorithm for computing edge lengths (Theorem \ref{thm:compute}) reduces, in the case of OLS, to that of \cite{Desper2002}. It has the same optimal running time as the algorithms in \cite{Bryant1998,Gascuel1997,Vach1989}.

\begin{example}[Balanced minimum evolution]
\end{example}
The {\bf B}alanced {\bf M}inimum {\bf E}volution model was introduced by Pauplin in \cite{Pauplin2000}. The motivation was that in the computation of $\hat{l}_{e^*}$ in the OLS model, the distances $D_{ac}$ and $D_{bd}$ can receive different weights than $D_{ad}$ and $D_{bc}$ where $a \in A, b \in B, c \in C$ and $d \in D$. Pauplin therefore suggested an alternative model where all clades are weighted equally. 
\begin{corollary}[Pauplin's edge formula]
The WLS edge lengths with variance model $V_{ij} \propto 2^{|\overline{i,j}|}$ are given by
$\hat{l}_{e^*}=\frac{1}{4}(D_{AC}+D_{BD}+D_{AD}+D_{BC})-\frac{1}{2}(D_{AB}-D_{CD})$
 for internal edges and $\hat{l}_{e^*}=\frac{1}{2}(D_{Ai}+D_{Bi})-\frac{1}{2}(D_{AB})$ for edges adjacent to leaves.
\end{corollary}
{\bf Proof:} This corresponds to the multiplicative variance model with $w_e=0.5$ for all edges $e$. One can easily show that in this case $Z_{AB}\propto 2^{-|\overline{A,B}|}$ and the result follows trivially from Theorem \ref{thm:main}. 

As far as we are aware, this is the first proof that the formulas given by Pauplin for edge lengths are in fact the WLS edge weights under the variance model described above. This implies:
\begin{remark}
\label{remark:nj_edgeweights}
The edge weights of the neighbor-joining tree obtained from the standard reduction formula are equal to the weighted least squares edge length estimates under the BME model.
\end{remark}
This result is a companion to the 
the connection between Pauplin's tree length formula and WLS tree length under the BME model that was established by Desper and Gasquel in \cite{Desper2004}. They proved the following:
\begin{corollary}[Desper and Gascuel \cite{Desper2004}]
The tree length estimator given by $\hat{l}=\sum_{ab}D_{ab}2^{1-p_{ij}}$ is the minimum variance tree length estimator for the BME model. It is also identical to the one given by the coefficients $p^t = f^t  (S_T^t  V^{-1}  S_T)^{-1}  S_T^t  V^{-1}.$
\end{corollary}
{\bf Proof:} The second part of the corollary follows trivially from Theorem \ref{thm:Gauss}. The first part follows from a simple combinatorial argument by adding up the WLS edge lengths.
Alternatively, one can notice directly that since $p_{ab}=2^{1-p_{ij}}$, it follows that $p_{ab}V_{ab}$ is the uniform vector, and thus defines a $T$-additive map, corresponding to the star topology (equal-length leaf edges and zero-length internal edges). Finally, $\sum_{i,j}S_{ij,e} p=1$ follows from an easy counting argument.
Further elaboration on Remark \ref{remark:nj_edgeweights} is beyond the scope of this paper.

\begin{example}[The taxon-weighted variance model]
\end{example}
Another well known WLS model was introduced by Denis and Gascuel in \cite{Denis2003}. Under this model we set $V_{ij}=t_i t_j$ for some $t_1,\ldots,t_n \in \mathbf{R}_+$. In the tree-multiplicative model, this corresponds to setting $w_e=1$ for internal edges and $w_e=t_i$ when $e$ is the leaf edge adjacent to leaf $i$. The paper \cite{Denis2003} gives a beautiful proof for the statistical consistency of this model (which implies statistical consistency of OLS), and also provides an $O(n^2)$ algorithm for computing the WLS edge lengths. However, the algorithm is based on a recursive agglomeration scheme and an explicit formula for the edge lengths based on the values of $D$ is not given. Such a formula follows from Theorem \ref{thm:main}:
\begin{corollary}
For $e$ an internal edge of $T$, the WLS edge length $\hat{l}_{e^*}$ is given by
\begin{eqnarray}\label{thm:WLSInternal2}
2\hat{l}_{e^*} & = & \frac{T_AT_D+T_CT_B}{(T_A+T_B)(T_CT_D)}(D_{AC}+D_{AC})\nonumber\\
& + & \frac{T_AT_C+T_DT_B}{(T_A+T_B)(T_CT_D)}(D_{AD}+D_{BC})\nonumber\\
& - &(D_{AB}+D_{CD})
\end{eqnarray}
where $T_X=\sum_{x\in X} t_x$ and $D_{XY}=\sum_{x\in X, y\in Y}\frac{t_xt_y}{T_XT_Y}D_{xy}$. If $e^*$ is adjacent to a leaf, 
\begin{equation}\label{thm:WLSLeaf2}
2\hat{l}_{e^*}=D_{Ai}+D_{Bi}-D_{AB}.
\end{equation}
\end{corollary}

\section{Final remarks}

An important question  is whether the variance matrices required for the IIP property to hold are realistic for problems where branch lengths are estimated using standard evolutionary models. In fact, semi-multiplicative matrices do not exactly capture the desired form of the variance, but they are good approximations. We illustrate this  for the Jukes--Cantor model \cite{Jukes1969}:
\begin{proposition}[Variance of distance estimates \cite{Bulmer1991,Nei1989}]
Let the random variable $Y$ be the fraction of different nucleotides between two sequences of length $n$ that are generated from the Jukes--Cantor process with branch length $\delta$. Then the expected value of the empirical distance $D=-\frac{3}{4}\log \left(1-\frac{4}{3}Y \right)$ is $\delta$ and its variance is
\begin{equation}
\label{eq:JCC}
Var(D) \approx \frac{3}{16n} \left( 3e^{\frac{8}{3}\delta} + 2e^{\frac{4}{3}\delta} -3 \right).
\end{equation}
\end{proposition}
This result can be extended to more general models. Since the branch lengths for an evolutionary model are tree-additive, this shows that for many regimes of the parameter $\delta$, a tree-multiplicative model for variances is very reasonable. For a discussion on the statistics rationale behind least squares see \cite{Felsenstein2003}.

Unfortunately, the Fitch--Margoliash assumption that the variance $V_{ij}=\mbox{Var}(D_{ij}) \propto D_{ij}^2$ is inaccurate in light of (\ref{eq:JCC}), nor does it lead to IIP estimates since $V$ is not semi-multiplicative. This means that for generic dissimilarity maps, the Fitch--Margoliash least squares estimates of edge lengths will depend on irrelevant distance estimates.

Another point that is important  is that although it follows
from Theorem \ref{thm:Gauss} that for any $V$ and $f$ there
is a unique BLUE $p$ for $f^t  l$, the converse of this
statement is not true. For example, if $p$ is BLUE for $f^t  l$
with variance matrix $V$, then $p$ is BLUE for $f^t  l$ with
variance matrix $kV$
where $k \geq 0$. This is obvious because $S_T^tp$ remains the same, and
$kV  p$ is a $T$-additive map if $V  p$ is a $T$-additive map. However this point has
more subtle (and serious) consequences:

\begin{proposition}[Non-uniqueness of tree length]
\label{prop:paupuniq}
The WLS estimated tree length with $V=(c_1+c_2(|\overline{i,j}|-1))2^{|\overline{i,j}|}$ does not depend on the constants $c_1$ and $c_2$.
\end{proposition}

Proposition \ref{prop:paupuniq} has significance for the interpretation of
the neighbor-joining algorithm. Based on \cite{Desper2004}, in \cite{Steel2006}  it is shown that
neighbor-joining minimizes the balanced evolution criterion at
each step. The criterion is argued to be statistically relevant by
virtue of the fact that it is the BLUE for the tree length under
the assumption that $V_{ij} \propto 2^{|\overline{i,j}|}$. Proposition \ref{prop:paupuniq}
shows that there are many (significantly) different variance assumptions
that yield the same tree length estimate. In fact, for some tree topologies, it is even possible that the OLS tree length is equal to the BME WLS tree length (for example for 5 taxa trees).
This means that by minimizing the tree length some
information about the variance is being discarded, and from this point
of view the fact that the balanced minimum evolution criterion is equal
to the BLUE tree length for multiple variance assumptions can be seen as a
weakness of balanced minimum evolution methods, not a strength.

There are other issues that are important in least squares applications in phylogenetics that we have not mentioned in this paper. One obvious difficulty with applying WLS methods to tree length estimation is that the resulting estimators are tree-additive, and not necessarily tree-metrics. That is, there may be edge length estimates that are negative. A number of strategies for solving the non-negative WLS problem have been proposed \cite{Felsenstein1997,Gascuel1996,Hubert1995, Makarenov1999}.

Our optimal algorithm for weighted least squares edge length estimates for multiplicative matrices is similar in spirit to a some of the algorithms in \cite{Bryant1998}. In fact, we believe that all the fast algorithms for WLS edge lengths can be understood within a single framework. The unifying concept is the observation that they all essentially estimate the Lagrange tree, either via a top-down, or bottom-up approach. We defer a detailed discussion of this to another paper. Finally, a key issue is that of consistency for specific forms of variance matrices assigned to all trees \cite{Denis2003,Willson2005}. An obvious question is what classes of semi-multiplicative variance matrices result in consistent tree estimates. A full discussion of this topic is also beyond the scope of this paper.
\section{Acknowledgments}
Radu Mihaescu was supported by a National Science Foundation Graduate Fellowship and partially by the Fannie and John Hertz foundation. Lior Pachter was supported in part by NSF grant CCF-0347992.



\bibliography{phylo}

\bibliographystyle{plain}











\end{document}